\font\smc=cmcsc10

\def \varpi {\bar \omega}
\def \N {I \! \! N}
\def \Z{Z \! \! \! Z}

\def \R {I \! \! R}
\def \N {I \! \! N}

\def\text#1{{\textstyle\rm#1}}

\def\newline{\unskip\null\hfill\break}
\def\square{\sqcup\!\!\!\!\sqcap}

\def\cejour{
\number\day\space\ifcase\month\or janvier\or f\'evrier \or
mars \or avril\or mai \or juin\or juillet\or ao\^ut\or
septembre\or octobre\or novembre\or
d\'ecembre\fi\space\number\year} 
\def\today{
\space\ifcase\month\or January\or February \or
March \or April\or May \or June\or July\or August\or
September\or October\or November\or
December\fi\space
\number\day ,\space\number\year} \nonstopmode
\def\cejour{(\number\time)
\number\day\space\ifcase\month\or janvier\or f\'evrier \or
mars \or avril\or mai \or juin\or juillet\or ao\^ut\or
septembre\or octobre\or novembre\or
d\'ecembre\fi\space\number\year} 
\def\C{{\mathchoice {\setbox0=\hbox{$\displaystyle\rm C$}\hbox{\hbox
to0pt{\kern0.4\wd0\vrule height0.9\ht0\hss}\box0}}
{\setbox0=\hbox{$\textstyle\rm C$}\hbox{\hbox
to0pt{\kern0.4\wd0\vrule height0.9\ht0\hss}\box0}}
{\setbox0=\hbox{$\scriptstyle\rm C$}\hbox{\hbox
to0pt{\kern0.4\wd0\vrule height0.9\ht0\hss}\box0}}
{\setbox0=\hbox{$\scriptscriptstyle\rm C$}\hbox{\hbox
to0pt{\kern0.4\wd0\vrule height0.9\ht0\hss}\box0}}}}
\def\cp#1{\bbbc{\rm P}^{#1}}

\def\frac#1#2{{#1\over#2}}
\def\square{\sqcup\!\!\!\!\sqcap}

\def\norm{\mid\!\mid}

\def\cp#1{\C{\rm P}^{#1}}


\noindent \today\vskip 1mm

\centerline {\bf Complex Plateau problem: old and new results and
prospects}\vskip 2mm

\centerline {Pierre Dolbeault}\vskip 3mm

\noindent {\it Abstract.} The Plateau problem is the research of a surface of minimal area, 
in the 3-dimensional Euclidean space, whose boundary is a given continuous closed curve. The complex Plateau
problem is analogous in a Hermitian complex manifold: it  is a geometrical problem of extension of a closed
real curve or manifold into a complex analytic subvariety, or into a Levi-flat subvariety. 
Wirtinger's inequality in $\C^n$ is recalled. Minimatity of complex analytic subvarieties and analogous properties of Levi-flat subvarieties, in K\"ahler manifolds, are
given. Known results in $\C^n$ and $\C P^n$ are recalled. 
Extensions to real parametric problems are solved or proposed, leading
 to the construction of Levi-flat hypersurfaces with prescribed
boundary in some complex manifolds.\noindent \vskip 2mm

\centerline {Contents}\vskip 1mm

1. Introduction

2. Volume minimality of complex analytic subvarieties and  of Levi-flat hypersurfaces in K\"ahler manifolds

3. Possible origin: holomorphic extension; polynomial envelope of a real curve

4. Solutions of the complex Plateau problem (or boundary problem) in different spaces

5. Extension to real parametric problems

6. Levi-flat hypersurfaces with prescribed boundary: preliminaries

7. Levi-flat hypersurfaces with prescribed boundary: particular cases\vskip 2mm

\noindent {\bf 1. Introduction.}\vskip 2mm

Given a Hermitian manifold $X$, the {\it complex Plateau problem} is the research of an even dimensional
subvariety with negligible singularities, with given boundary, and of minimal volume in $X$.
We will call {\it mixed Plateau problem} the research of a real hypersurface with given
boundary, and of minimal volume in
$X$. More briefly, both problems will be called complex Plateau problem. 

First we shall recall or show that complex analytic subvarieties, resp. Levi-flat hypersurfaces are solutions
of the Plateau problem when $X$ is K\"ahler (section 2). 

Then we will consider the complex Plateau problem as the
research of the extension of an odd dimensional, compact, oriented, connected submanifold into a
complex analytic subvariety, and recall known solutions (section 3, 4). 

To solve the mixed Plateau problem as the research of the
extension of an oriented, compact, connected, 2-codimensional submanifold into a Levi-flat 
hypersurface, we
will need solutions of the complex Plateau problem with real parameter, in $\C^n$ and $\cp n$; in
$\cp n$, it is an open problem to explicit satisfactory conditions on the boundary. In this way, we get 
very peculiar solutions of mixed Plateau problems (section 5).

Finally, the mixed Plateau problem is solved in $\C^n$, in particular cases, as  a projection of a Levi-flat variety, and set up in $\cp n$ (sections
6,7): known solutions are recalled in $\C^n$ when the complex points of 
the boundary are elliptic; special elliptic and hyperbolic points of the
boundary are defined, and a solution when the boundary is a "horned
sphere" is described; this will be the opportunity to precise and complete
results announced in ([D 08], section 4). Problems when the boundary has
general hyperbolic points are still open. 

Proofs of the results in sections 6 and 7 will appear in detail elsewhere
[D 09].

\noindent {\bf Acknowledgments}. I thank G. Tomassini and D. Zaitsev for 
discussions, corrections and remarks about several parts of this paper .

\vskip 2mm

\noindent {\bf 2. Volume minimality of complex analytic subvarieties and  of Levi-flat hypersurfaces in K\"ahler manifolds.} 

\vskip 1mm

\noindent 2.1. {\it Wirtinger's inequality} (1936)  [H 77].

In $\C^n$,with complex coordinates $(z_1,\ldots, z_n)$, we have the Hermitian metric $H=\displaystyle\sum_{j=1}^n dz_j\otimes d\overline z_j$
and the exterior form (standard K\"ahler) $\omega=\displaystyle\frac{i}{2}\sum_{j=1}^n dz_j\wedge
 d\overline
z_j=\displaystyle\sum_{j=1}^n dx_j\wedge dy_j$.

From the real vector space $\R^{2n}\cong\C^n$, we consider the real vector
 space $\Lambda_{2p}\R^{2n}$ of the $2p$-vectors with the associated norm
$\vert .\vert$; every decomposable vector (exterior product of elements
of $\R^{2n}$) defines a real $2p$-plane of
$\C^n$ i.e. an element of the Grassmannian $G_{2n}^{2p}$. We define the
norm
$\norm\zeta\norm=\displaystyle
\inf\sum_j\mid\zeta_j\mid$ where $\zeta=\displaystyle\sum_j\zeta_j$,
$\zeta_j$ is decomposable.

Let $P_{pp}=\{\displaystyle\sum_{j=1}^N\lambda_j\zeta_j; \zeta_j$ decomposable defining a complex $p$-plane of
$\C^n$; $\lambda_j\geq 0;N\in \N^*\}$.

\noindent 2.1.1. {\bf Theorem.} {\it For every $\zeta\in
\Lambda_{2p}\C^n$, we have:

$$\frac{1}{p!}\omega^p(\zeta)\leq\norm\zeta\norm;$$
equality uniquely for  $\zeta\in P_{pp}$} $\lbrack$W 36$\rbrack$.\vskip
1mm

\noindent 2.1.2. {\bf Corollary.} {\it Let $V$ be a smooth real oriented
$2p$-dimensional submanifold of a Hermitian manifold $X=(X,\omega)$ of
complex dimension $n$. Then $\int_V\omega^p/p!\leq {\rm vol}_{2p}(V)$
with equality iff $V$ is complex.}\vskip 1mm

\noindent 2.2. {\it
Currents with measure coefficients.} [H 77]     \vskip 1mm

\noindent 2.2.1. {\it Comass of an $r$-form; mass of a
current with measure coefficients.}  

Let $\varphi\in
\Lambda^r\R^{2n}$, the comass of $\varphi$ is defined as
$$\norm\varphi\norm^*={\rm sup}\{\varphi(\zeta): \zeta\in
G^r_{2n}\subset\Lambda_{2p}\R^{2n}\}$$

Let $\Omega$ be an open subset of $\C^n$, for every differential form
$\varphi$ of degree $r$ on $\Omega$, let 
$$\norm\varphi\norm^*={\rm sup}\{\norm\varphi(z)\norm^*: z\in \Omega\}$$
where $\norm\varphi(z)\norm^*$ is the comass of $\varphi (z)$.

Let $T$ be a  current with measure coefficients on $\Omega$, $K$ be
any compact subset of $\Omega$ and $\chi_K$ the characteristic function
of $K$,
$$M_K(T)= \sup_{\norm\varphi\norm^* \leq 1}\vert\chi_KT(\varphi)\vert$$
is, by definition,  the {\it mass} of $T$ on $K$.\vskip 1mm

The measure which assigns the number $M_K(T)$ to each compact set $K\subset\Omega$ is called the {\it mass} or {\it volume measure} of $T$ and denoted $\norm T\norm$, so that $M_K(T)=\norm T\norm (K)$.\vskip 1mm

 \noindent 2.3. {\it Complex Plateau problem}.\vskip 1mm

\noindent 2.3.1. [H 77] On $\Omega\subset\C^n$, or more generally, on a 
Hermitian manifold $(X,\omega)$,  let $B$ be a $d$-closed current of
dimension $2p-1$ with compact support, and let $T$ be a $(2p)$-current
with compact support and measure coefficients such that $dT=B$.      The
complex Plateau problem is to find such a $T$ with minimal mass,  i.e.
for every compactly supported current $S$, with measure coefficients such that $dS=B$, to have
$M(T)\leq M(S)$, or equivalently, for every compactly supported,
$d$-closed $(2p)$-current with measure coefficients $R$, 
$$M(T)\leq M(T+R)$$ 
Such a $T$ is said {\it absolutely volume minimizing} on $X$.\vskip 1mm

Let $T$ be a $d$-closed $(2p)$-current with measure coefficients on $X$.
If, for each compact subset $K$ of $X$,
$$M_K(T)\leq M(\chi_KT+R)$$for all compactly supported $d$-closed
$(2p)$-current $R$ with measure coefficients on $X$, then $T$ is said to 
be {\it absolutely volume minimizing} on $X$.

\vskip 1mm
\noindent 2.3.2. {\bf Theorem.} [H 77] {\it Let $T$ be a $2p$-current
with measure coefficients on a Hermitian manifold $(X, \omega)$ and $K$
be a compact subset of $X$. Then
$$(\chi_KT)(\omega^p/p!)\leq M_K(T).$$
and equality holds iff $\chi_KT$ is strongly positive.}\hskip 10mm
$\square$\vskip 1mm

\noindent 2.3.3. [H 77] {\it Volume minimality of complex analytic sets 
in a K\"ahler manifold.}\vskip 1mm

\noindent 2.3.4. {\bf Corollary to Theorem 2.3.2.} {\it Assume that
$X=(X,\omega)$ is a K\"akler manifold and does not contain compact
$p$-dimensional complex subvarieties. Let $V$ be a $p$-dimensional
complex subvatiety, and
$T=\lbrack V\rbrack$, then $T$ is absolutely volume minimizing} on
$X$.     

\noindent {\it Proof}.  $T$ is strongly positive. Let $K$ be a compact
subset of $X$ and $R$ be a compactly supported $d$-closed
$(2p)$-current with measure coefficients. From Theorem 2.3.2, $M_K(T)=(\chi_KT)(\omega^p/p!)
$. But locally $\omega=dd^c\psi$, then $\omega^p=\omega^{p-1}\wedge
dd^c\psi=d(\omega^{p-1}\wedge d^c\psi)$, so in the neighborhood of any point
of $X$, $R(\omega^p)= R(d(\omega^{p-1}\wedge d^c\psi))$.
Let $(\alpha_j)_{j\in J}$ be a partition $C^\infty$ of unity subordinate to a locally finite open covering $(U_j)_{j\in J}$ of $X$ such that for every $j$, $\omega\vert_{U_j}=dd^c\psi_j$. Then 
 $$R(\omega^p)=\sum_j\alpha_j R(d(\omega^{p-1}\wedge d^c\psi_j))=\pm\sum_jd(\alpha_j R)(\omega^{p-1}\wedge d^c\psi_j)=0,$$
 because: $$\sum_jd(\alpha_j R)=\sum_jd\alpha_j\wedge R + \sum_j\alpha_j\wedge dR=0$$
 and, as in the proof of ([H 77], Corollary 1.25), in an open set of the Hermitian $\C^n$, 
$$M_K(T)= (\chi_KT)(\omega^p/p!)=(\chi_KT+R)(\omega^p/p!)
\leq M_K(T+R).\eqno\square$$\vskip 1mm

\noindent 2.3.5. {\bf Remark.} If $X$ contains a compact $p$-dimensional complex
subvariety $W$, $d\lbrack V\rbrack=0$, but $M_K(\lbrack V\rbrack >0$; then $T$ is {\it relatively volume minimizing} on $X$.     \vskip 2mm

\noindent 2.4. {\it Volume minimality of Levi-flat hypersurfaces in K\"ahler manifolds.}

   We suppose to be in the category of currents with measure coefficients.\vskip 1mm
   
   Recall the definition: A {\it Levi-flat subvariety} (with
negligible  singularities), of odd dimension, is, outside of the
singularities, a submanifold with Levi form $\equiv 0$, or, equivalently, is
foliated by complex analytic hypersurfaces.

Let $M$ be a $C^\infty$ Levi-flat hypersurface of a $C^\infty$  K\"ahler manifold $X=(X,\omega)$ bearing a foliation ${\cal L}$ by complex hypersurfaces $M_l$ and let $L$ be the space of the foliation ${\cal L}$ assumed to be a $C^\infty$ real curve.

Let $M'$ be a $C^\infty$ hypersurface of $X$ bearing a foliation ${\cal L}'$ with the same space $L$; the leaves of ${\cal L}'$ being $C^\infty$ subvarieties with negligible singularities.

Let $S$ be a $C^\infty$ compact submanifold of codimension 2 of $X$. We denote by the same notation the hypersurfaces and submanifold and the integration currents they define.
\vskip 1mm

\noindent 2.5. {\it Mixed Plateau problem}.\vskip 1mm

Given $S$ to find a $C^\infty$ hypersurface in $X\setminus S$ whose boundary is $S$ in the 
category ${\cal H}$ of foliated hypersurfaces with the same space of foliation, a real
curve. If $M'$ is such a hypersurface whose space of foliation is $L$ and the leaves
$(M'_l, l\in L)$, then vol$(M')=\int_L{\rm vol}(M'_l)dl$.

From section 2, for every $l\in L$, vol$(M'_l)\geq {\rm vol}(M_l)$ then {\it $M$ is relatively volume minimizing in
the category ${\cal H}$} and, by definition,
 {\it $M$ is solution of the mixed Plateau problem.}\vskip 1mm

\noindent 2.6. {\it Research of solutions of the complex Plateau
problem}.\vskip 1mm

The present method of resolution consists in finding complex analytic,
resp. Levi-flat subvarieties, in $X\setminus S$, whose boundary $S$ (in
the sense of currents) is a submanifold of $X$ with convenient
properties. .
\vskip 2mm

\noindent {\bf 3. Possible origin: holomorphic extension; polynomial envelope of a real curve.}\vskip 1mm

 \noindent 3.1. The extension theorem of Hartogs, obtained at the beginning of the $20^{th}$ century, has been completely proved by Bochner and  Martinelli, independently, in 1943. 
 The simplest version is:
 
 {\it Let $\Omega$ be a bounded open set of  $\C^n$, $n\geq 2$. Suppose that $\partial\Omega$ be of class $C^k$ $(1\leq k\leq\infty)$ or of class $C^\omega$ (i.e. real analytic). Let  $f$ be a function in $C^l(\partial\Omega)$, $1\leq l\leq k$. 

Then the two conditions are equivalent:

$(i)$ $f$ is a CR function, i.e. the differential of $f$ restricted to the complex subspaces of the  
tangent space to 
 $\partial\Omega$, at every point, is $\C$-linear;

$(ii)$ there exists $F\in C^l(\overline\Omega)\cap {\cal O}(\Omega)$ such that
$F\mid_{\partial\Omega}=f$}.\vskip 1mm
 
 Then the graph of $f$ is the boundary of the 
 complex analytic submanifold defined by the graph of $F$ in $\C^{n+1}$.\vskip 1mm
 
 \noindent 3.2. Let $M$ be a compact submanifold of  dimension 1 of $\C^n$, we call
polynomial envelope of $M$, the compact set $\{z\in\C^n; \displaystyle\mid P(z)\mid\leq \max_{\zeta\in M}\mid 
P(\zeta)\mid ;
P\in\C\lbrack z\rbrack$, the polynomial ring with complex coefficients $\}$. 
 
 Then (J. Wermer
(1958)), the polynomial envelope of
$M$ is either $M$, or the union of $M$ with the support of a complex analytic variety $T$, of
complex dimension 1, whose boundary is $M$ $\lbrack$We 58$\rbrack$.

 \vskip 2mm

\noindent {\bf 4. Solutions of the complex Plateau problem (or boundary problem) in different spaces.}\vskip 2mm

\noindent 4.1. The first result has been obtained in 1958, by J. Wermer, 
in $\C^n$, for $p=1$ and $M$ holomorphic image of the unit circle in $\C$
$\lbrack$We 58$\rbrack$; this result has been generalized to the case where
$M$ is a union of $C^1$ real connected curves by Bishop, Stolzenberg (1966),
looking for the polynomial envelope of $M$ according to section 3.2.
\vskip 1mm

In $\C^n$, after preliminary results by Rothstein (1959) $\lbrack$Rs 59$\rbrack$
, the boundary problem has been solved by Harvey and Lawson (1975), for $p\geq
2$, under the necessary and sufficient condition: $M$ is compact, maximally
complex and, for $p=1$, under the moment condition: $\int_M\varphi=0$,  for
every holomorphic 1-form $\varphi$ on $\C ^n$ $\lbrack$ HL $75 \rbrack$. For
$n=p+1$, the method, inspired by the  Hartogs' theorem consists in building
$T$ as the divisor of a meromorphic function the {\it defining function} $R$;
this function itself is  constructed, step by step, from solutions of
$\overline\partial$-problems with compact support. $T$ can also be viewed as
graph (with multiplicities on the irreducible components) of an analytic
function with a finite number of determinations. 
 For any $p$, we come back to the particular case using  projections.\vskip 1mm
 
 In $\cp n\setminus\cp {n-r}$, $1\leq r\leq n$, for compact $M$, the problem 
has a une solution if and only if, for $p\geq r+1$, $M$ is maximally complex
and if, for $p=r$, $M$ satisfies the moment condition:  $\int_M\varphi=0$, for
every $\overline\partial$-closed $(p,p-1)$-forme $\varphi$
. The method consists in solving the boundary problem, in $\C^{n+1}\setminus\C
^{n-r+1}$, for the inverse image of $M$ by the canonical projection 
$\lbrack$HL $77 \rbrack$.

In both cases, the solution is unique.

Harvey et Lawson assume the given $M$ to be, except for a closed set of Hausdorff $(2p-1)$-dimensional measure zero, an oriented manifold of class $C^1$; we will say: $M$ is a {\it variety $C^1$ with negligible singularities}.

The boundary problem in $\cp n$ has been set up, for the first time, by J. King $\lbrack$ Ki $79 \rbrack$; uniqueness of the solution is no more possible, since two solutions differ by an algebraic $p$-chain.

4.2. In $\cp n$, a solution of the boundary problem has been obtained by P. Dolbeault et G. Henkin 
for $p=1$, (1994), then for every $p$ (1997) and more
generally, in a $q$-linearly concave domain $X$ of $\cp n$, i.e. a union of projective subspaces of dimension $q$ $\lbrack$DH 97 $\rbrack$.

The necessary and sufficient condition for the existence of a  solution is an  extension of the moment condition: it uses a Cauchy residue formula in one variable and a non linear differential condition which 
appears in many questions of Geometry or mathematical Physics. In the simplest case: $p=1$, $n=2$, this is the shock wave equation for a local holomorphic function in 2 variables $\xi, \eta$, \hskip 3mm $\displaystyle f{\frac{\partial f}{\partial\xi}}={\frac{\partial f} {\partial\eta}}$.

 From a local condition, the above relation allows to construct, by extension 
ot the coefficients, a meromorphic function playing, in $\C ^n$, the same part
as the  Harvey-Lawson defining function described above; it defines a
holomorphic
$p$-chain extendable to
$\cp n$ using the classical Bishop-Stoll theorem.

4.2.1. The conditions of regularity of $M$ have been weakened, first in $\C ^n$, and 
for $p=1$, to a condition, a little stronger than the rectifiability, by H.
Alexander
$\lbrack$Al 88 $\rbrack$ who, moreover, has given an essential counter-example $\lbrack$ Al 87$\rbrack$, then by Lawrence $\lbrack$Lce 95$\rbrack$ and finally, and for any $p$, in
$\C ^n$ and $\cp n$, by T.C. Dinh $\lbrack$Di 98$\rbrack$: $M$ is a
rectifiable current whose tangent cone is a vector subspace almost everywhere.
Moreover, Dinh has obtained the reduction of the boundary problem in $\cp n$
to the case $p=1$, with weaker conditions than above and by an elementary
analytic procedure $\lbrack$Di 98$\rbrack$.
 
 All the previous results are obtained as  Corollaries.
 
 New progress by Harvey and Lawson [HL 04].\vskip 2mm
 
 \noindent {\bf 5. Extension to real parametric problems.}\vskip 2mm

 \noindent 5.1. {\it In a real hyperplane of $\C^n$.}\vskip 1mm
 
 5.1.1. Let $E\cong{\bf R}\times\C^{n-1}$, and $k:{\bf
R}\times\C^{n-1}\rightarrow {\bf R}$ be the projection. 
Let $N\subset E$ be a
compact, (oriented)  CR subvariety of $\C^n$ of
real dimension $2n-4$ and CR dimension
$n-3$, $(n\geq 4)$, of class $C^\infty$, with negligible
singularities (i.e. there exists a closed subset $\tau\subset N$ of
$(2n-4)$-dimensional Hausdorff measure $0$
such that $N\setminus \tau$ is a CR submanifold).
Let $\tau'$ be the set of all points $z\in N$
such that either $z\in\tau$ or $z\in N\setminus\tau$
and $N$ is not transversal to the complex hyperplane 
$k^{-1}(k(z))$ at $z$.
Assume that $N$,  as a
current of integration, is $d$-closed and satisfies:

{\rm (H)} there exists a closed subset $L\subset{\bf R}_{x_1}$ 
with
$ H^1(L)=0$ such that for every $x\in k(N)\setminus L$,
the fiber $k^{-1}(x)\cap N$ is connected 
and does not intersect $\tau'$.\vskip 1mm

\noindent  5.1.2. {\bf Theorem} [DTZ 09] (see also [DTZ 05]). {\it Let $N$
satisfy {\rm (H)} with
$L$ chosen accordingly. Then, there exists, in $E'= E\setminus k^{-1}(L)$, a
unique
$ C^\infty$ Levi-flat $(2n-3)$-subvariety $M$ with negligible singularities in
$E'\setminus N$, foliated by complex $(n-2)$-subvarieties, with the properties
that
 $M$ simply (or trivially) extends  
to
$E'$ as a $(2n-3)$-current  (still denoted $M$) such that $dM=N$ in  
$E'$.
The leaves are the sections by the hyperplanes $E_{x_1^0}$, $x_1^0\in
k(N)\setminus L$, and are the solutions of the ``Harvey-Lawson problem'' for  
finding a holomorphic subvariety in $E_{x_1^0}\cong\C^{n-1}$ with prescribed boundary $N\cap E_{x_1^0}$.}\vskip 1mm

\noindent  5.2. {\it In a real hyperplane of $\cp {n+1}$.}\vskip 1mm

\noindent  5.2.1. 
  The simplest significant case is the boundary problem in $\C P^3$. For the boundary problem with real parameter in $\C^3$, we considered a boundary problem in $\R\times\C^3$, i.e. in the subspace of $\C^4$, in which the first coordinate is real. In the same way, we will consider in $\C P^4$, with homogeneous coordinates $(w_0,w_1,\ldots,w_4)$, a boundary problem in the subspace $E$ defined by $w_1=\lambda w_0$, with $\lambda\in\R$. Then,  for personal convenience, we will follow, step by step, the known construction in $\C P^3$ in the oldest version $\lbrack$DH 97$\rbrack$.
  
  Particularly, the coefficients $C_m$ of the defining function $R$ of the solution are estimated as for the problem in $\C P^3$.
   
   The end of the proof of the main theorem seems analogous to the known case in $\R\times\C^3$.\vskip 1mm
  
\noindent  5.2.2. The projective space $\C P^3 $ has homogeneous  coordinates $w'=(w_0,w_2,\ldots,w_4)$; denote $Q=\{w_0=0\}$
the hyperplane at infinity of $\C P^3$. 
\vskip 1mm

For $w_0\not =0$, let $k$ be the projection: $E\rightarrow \R_\lambda$, 
$(w_0,w_1=\lambda w_0, w_2,w_3,w_4)\mapsto\lambda$; for $w_0=0$, $\lambda$ is 
indeterminate..
We also have the projection: \hskip 5mm$\pi: E\rightarrow\C P^3$,
$(w_0,\lambda w_0, w_2,w_3,w_4)\mapsto (w_0,
w_2,w_3,w_4)$. In the same way,
$(E\setminus\{w_0=0\}) \cong\R\times\C^3 $.

\noindent 5.2.3. Let $N\subset E\subset\C P^4$ be a submanifold of class
$C^\infty$, CR,
oriented, compact of $E$, of dimension 4,  of CR dimension 1, with negligible 
singularities. $N$ being compact in $E$, $k(N)$ is compact in $\R$, i.e.  in $N$, the
parameter $\lambda$ varies in a closed, bounded interval $\Lambda$ of $\R$. 

Assume that $N$ satisfies the same properties as in subsection 5.1.1.\vskip 1mm

\noindent 5.2.4. Consider the complex hyperplanes of $\C P^4$, whose equation is
$\tilde h(w)=w_4-\xi_2 w_0-\eta'_2 w_1-\eta_2 w_2=0$
and, in $E$, the subspaces $P_{\nu'_\lambda}$ whose equation is 
  $$\tilde h_1(w',\lambda)=w_4-\xi_2 w_0-\eta'_2\lambda w_0-\eta_2
w_2=w_4-(\xi_2 +\eta'_2\lambda) w_0-\eta_2 w_2=0,$$ of real dimension
 5. Their restrictions to $(E\setminus\R\times Q)
\cong\R\times\C^3 $ are real affine sub-spaces of dimension 5.

We note $\nu'_\lambda$ the $1\times 2$-matrix
$\bigl((\xi_2+\eta'_2\lambda)\  \eta_2\bigr)$.
  
  Generically, $\Gamma_{\nu'_\lambda}=N\cap
P_{\nu'_\lambda}$ is of dimension 2.
  
  For $z\in N, \lambda=k(z)$. let $E_\lambda=k^{-1}k(z)$; for $\lambda
\notin L$,
$N\cap E_\lambda$ is of dimension 3 and is contained in $E_\lambda\cong
\C P^3$.
  
  Consider the linear forms
$$\underline h(w)=w_3-\xi_1 w_0-\eta'_1 w_1-\eta_1 w_2$$
  
  $$\tilde h_0(w',\lambda)=w_3-\xi_1 w_0-\eta'_1\lambda w_0-\eta_1
w_2=w_3-(\xi_1+\eta'_1\lambda) w_0-\eta_1 w_2$$ 
Denote $\nu_\lambda=(\xi_\lambda,\eta)$ the $2\times 2$-matrix
$$\pmatrix{
(\xi_1+\eta'_1\lambda)  &\eta_1  \cr 
(\xi_2+\eta'_2\lambda) &\eta_2\cr
}$$
For fixed $\lambda$, $\nu_\lambda$ is a coordinate system of a chart of 
the Grassmannian $G_{\C}(2,4)$, i.e. $\nu_\lambda$ is a coordinate system of a
chart of  $G_{\C}(2,4)\times\R$. and we identify
$\nu_\lambda$ with the point of $G_{\C}(2,4)\times\R$ having these
coordinates.\vskip 1mm

Let  
$\xi_\lambda=^t\bigl((\xi_1+\eta'_1\lambda)\
(\xi_2+\eta'_2\lambda)\bigr)$;\hskip 3mm $\eta=^t(\eta_1\ \eta_2)$.
Remark that $\xi_\lambda$ depends on
$(\xi_1,\xi_2,\eta'_1,\eta'_2)$; to get
effective dependance on the  parameter $\lambda$, it suffices to fix
$\eta'_1\not =0$,
$.\eta'_2\not =0$. 

Recall:
$\xi_\lambda=^t(\xi_{\lambda 1}\ 
\xi_{\lambda 2}),$\hskip 3mm $\xi_{\lambda l}=\xi_l+\eta'_l \lambda$,
$l=1,2$, \hskip 3mm $\eta=^t(\eta_1\ \eta_2)$.
\vskip 1mm

   Let $z_j=w_j/w_0$, $j=2,\ldots, 4$, be the non
homogeneous coordinates ; $\tilde h_0$ defines the affine function:
      $$h= z_3-(\xi_1+\eta'_1\lambda)-\eta_1 z_2$$
      
     The two forms  $\tilde h_0$ et $\tilde h_1$ are linearly independent, then 
the set of their common zeros $D_{\nu_\lambda}$ is of real dimension 3, is
contained in
$P_{\nu'_\lambda}$; in general,
$D_{\nu_\lambda}\cap N$ is a finite set 
 $Z_{\nu_\lambda}$ ; then, for general enough fixed $\lambda$
and $\nu_\lambda$, 
$Z_{\nu_\lambda}=\emptyset$.
For every fixed real number $\lambda\notin L$, the situation in $E_\lambda$ is 
the classical situation in $\C P^3$.

\vskip 1mm
\noindent 5.2.5. {\it Boundary problem.}  Given $N$, find a complex 
analytic subvariety $M$ depending on the real parameter $\lambda$ such
that $dM=N$ in the sense of currents, under a necessary and sufficient
condition on $N$. 

To do this, we can check, step by step, the solution of the boundary
problem in $\C P^3$ $\lbrack$HL 97$\rbrack$, introducing the parameter
$\lambda$.

For $\lambda\notin L$, $\gamma_{\nu'_\lambda}=N\cap
P_{\nu'_\lambda}\cap E_\lambda$ is of dimension 1. Under the notations of the 
sub-section 5.2.4, consider the function 
$$G(\nu_\lambda)=\frac{1}{2\pi i}\int_{\gamma_{\nu'_\lambda}} z_2
\frac{dh}{h}\leqno (1)$$
\vskip 5mm

\noindent 5.2.6. {\bf Tentative statement.} {\it The following two conditions are 
equivalent:\vskip 1mm

$(i)$ There exists, in $E'=E\setminus k^{-1}L$, a $C^\infty$ Levi-flat 
subvariety $M$, (with negligible singularities), of dimension
$5$, foliated by complex analytic subvarieties $M_\lambda$ of complex  
dimension $2$, such that $M$
extends simply (or trivially) to $E'$ as a current of dimension
$5$
$($still denoted $M)$ such that $dM=N$ in $E'$. The leaves are the sections by 
the subspaces $E_\lambda, \lambda\in k(N)\setminus L$, and are the solutions 
of the boundary problem for finding complex analytic subvarieties in 
$E_\lambda\cong\C P^3$ with given boundary $N\cap
E_\lambda$.\vskip 1mm

$(ii)$ $N$ is a submanifold CR, oriented, of CR dimension $1$ outside a closed set of 4-dimensional Hausdorff measure $0$.

There exists a matrix $\nu_{\lambda^*}^*$ in the neighborhood of which
$$D^2_{\xi_\lambda}G(\nu_\lambda)=D^2_{\xi_\lambda}\sum_{j=1}^N
f_j(\nu_\lambda)$$
where $f_j$, $j=1,\ldots,N$, is a holomorphic function in $\nu_\lambda$,
$C^\infty$ en $\lambda$, and satisfies the system of P.D.E.  
$$f_j\frac{\partial f_j}{\partial\xi_{\lambda l}}=\frac{\partial
f_j}{\partial\eta_l}, \hskip 2mm l=2,3\leqno (2)$$}\vskip 1mm

\noindent 5.2.7. {\it Remark.} This result is not satisfactory because 
the 
relation of the analytic conditions with the geometry of the submanifold 
$N$ is not explicit.\vskip 1mm

\noindent 5.3. {\it Boundary problem in a real hyperplane of $\C^{n+1}$ or $\C
P^{n+1}$.}
\vskip 1mm

$\C^{n+1}$ and $\C P^{n+1}$ are both K\"ahler. {\it The solutions of the above
boundary problems are} both Levi flat, hence, from a plain extension of
section 2.5, volume minimal, i.e.
{\it solution, of codimension 3, of mixed Plateau problems}.\vskip 2mm

\noindent {\bf 6. Levi-flat hypersurfaces with prescribed boundary: preliminaries.}\vskip 2mm

\noindent 6.1. {\it Introduction.}

Let $S\subset \C^n$ be a compact connected 2-codimensional
submanifold. Find a Levi-flat
hypersurface $M\subset \C^n\setminus S$ such that $dM=S$
(i.e. whose boundary  is $S$, possibly as a current).

For $n=2$, near an elliptic complex point $p\in S$, $S\setminus
\{p\}$ is foliated by smooth compact real curves which bound analytic
discs (Bishop  [Bi 65]). The family of these discs fills  a smooth
Levi-flat hypersurface.

In 1983, Bedford-Gaveau considered the case of a particular
sphere with two elliptic
complex  points. If $S$ is contained in the boundary of a strictly
pseudoconvex bounded
domain, then the families of analytic discs in the neighborhood of
each elliptic point extend
to a global family filling a 3-dimensional  ball $M$ bounded by $S$.
In 1991, Bedford-Klingenberg [BeK 91] and Kruzhilin extended the result when there exist hyperbolic complex points on $S$ with
the same global condition.

Results of increasing 
generality have been obtained by Chirka, Shcherbina, Slodowski, G.
Tomassini until 1999. The global sufficient condition of embedding of
$S$ in the boundary of a strictly pseudoconvex domain is still required
in these papers. 

\vskip 1mm

 A first result for $n\geq 3$ (in the sense of currents),
and for elliptic points only, has been obtained four years ago ([DTZ 05]
and [DTZ 09] in detailed form); we got new results when
 $S$ is homeomorphic to a sphere, with three elliptic and
one hyperbolic special points (see [D 08] for a first draft), or
 a torus, with two elliptic and two hyperbolic special points and, more
generally, a manifold which is obtained by {\it gluing together 
elementary models}.

A local condition is required because,
in general, $S$ is not locally the boundary of a Levi-flat
hypersurface. The proof uses the construction of a foliation of $S$ by CR
orbits, Thurston's stability theorem for foliations on $S$, and a
parametric version of the
Harvey-Lawson theorem on boundaries of complex analytic varieties.
There is no global condition.\vskip 1mm

\noindent 6.2. {\it Preliminaries and definitions.}\vskip 2mm

\noindent 
6.2.1. A smooth, connected, CR submanifold $M\subset \C^n$ is called  {\it minimal} 
at a point
$p$ if there does not exist a submanifold
$N$ of $M$ of lower dimension through $p$ such that $HN = HM|_{N}$.
By a theorem  of Sussman,
all possible submanifolds $N$ such that $HN = HM|_{N}$ contain, at $p$, one of the
minimal possible dimension, called a CR {\it orbit} of $p$ in $M$
whose germ at $p$ is uniquely determined.\vskip 1mm

\noindent 6.2.2. $S$ is said to be a {\it locally flat boundary} at a
point
$p$ if it locally bounds a Levi-flat hypersurface near $p$. Assume that
$S$ is  CR in a small enough neighborhood $U$ of  $p\in S$. If all CR
orbits of
$S$ are $1$-codimensional (which will appear as a necessary condition
for our problem), the following two conditions are equivalent  [DTZ 05]:

$(i)$ $S$ is a locally flat boundary on $U$;

$(ii)$ $S$ is nowhere minimal on $U$.\vskip 1mm

\noindent 6.2.3. {\it Complex points of $S$} [DTZ 05]. 

At such a point
$p\in S$, $T_pS$ is a complex hyperplane in $T_p\C^n$. In suitable
holomorphic coordinates $(z,w)\in \C^{n-1}\times\C$ vanishing at $p$,
$S$ satisfies
$$ w= Q(z) + O(|z|^3),\quad
Q(z)= \sum_{1\leq i,j\leq n-1} (a_{ij}z_iz_j + b_{ij}z_i\overline z_j
+c_{ij}\overline z_i\overline z_j)\leqno (1)$$
$S$ is said {\it flat} at a
complex point $p\in S$ if \ $\sum b_{ij}z_i\overline z_j\in\lambda{\bf R},
\lambda\in {\C}$. We also say that $p$ is flat.

{\it Let $S\subset\C^n$ be a
locally flat boundary with a complex point $p$. Then $p$ is flat.}

By making the change of coordinates $(z,w)\mapsto(z,\lambda^{-1}
w)$, we make $\sum b_{ij}z_i{\overline z}_j\in\R$ for all $z$.
By a change of coordinates $(z,w)\mapsto (z,w+\sum
a'_{ij}z_iz_j)$ we can choose the holomorphic term in (1) to be
the conjugate of the antiholomorphic one and so make the whole form
$Q$ real-valued.

 We say that $S$ is in a {\it flat normal form} at $p$
if the coordinates $(z,w)$ as in (1) are chosen
such that $Q(z)\in{\bf R}$ for all $z\in\C^{n-1}$.\vskip 1mm

\noindent 6.2.4. {\it Properties of $Q$.}

Assume that $S$ is in a flat normal form; then, the quadratic form $Q$ is
real valued. Only holomorphic linear changes of coordinates are allowed.
If
$Q$ is positive definite or negative definite, the point $p\in S$ is said to be {\it
elliptic}; if the point $p\in S$ is is not elliptic, and if $Q$ is non degenerate, $p$ is said to be {\it hyperbolic}. From section 6.4, we will only consider particular cases of the quadratic form $Q$.\vskip
1mm

From [Bi 65], for $n=2$, in suitable holomorphic coordinates, \ \ 
$Q(z)=(z\overline z+\lambda Re \
z^2), \ \ \lambda\geq 0 $, under the notations of [BeK 91]; for
$0\leq\lambda<1$, $p$ is said to be {\it elliptic}, and for $1<\lambda$,
it is said to be {\it  hyperbolic}. The parabolic case $\lambda=1$, not
generic, is omitted [BeK 91]. When $n \geq 3$, the Bishop's result is not valid in general.\vskip 1mm

\noindent 6.3. {\it Elliptic points}.
\vskip 1mm

\noindent 6.3.2. {\smc Proposition} ([DTZ 05], [DTZ 09]).
{\it Assume that $S\subset\C^n$, ($n\ge 3$) is nowhere minimal at all
its CR points
and has an elliptic flat complex point $p$. Then there exists a
neighborhood $V$ of $p$ such that $V\setminus \{p\}$ is foliated by
compact real $(2n-3)$-dimensional CR orbits diffeomorphic to the
sphere ${\bf S}^{2n-3}$ and there exists a smooth function $\nu$,
 having the CR orbits as the level surfaces.}

\vskip 1mm

\noindent {\it Sketch of Proof} (see [DTZ 09]).

In the case of a quadric $S_0$ ($w=Q(z)$), the CR orbits are defined by
$w_0=Q(z)$, where $w_0$ is constant. Using (1), we approximate the
tangent space to $S$ by the tangent space to $S_0$ at a point with the
same coordinate $z$; the same is done for the tangent spaces to the CR
orbits on $S$ and $S_0$; then we construct the global CR orbit on $S$
through any given point close enough to $p$. 

\vskip 1mm

\noindent 6.4. {\it Special flat complex points}. We say that the flat complex
point
$p\in S$ is {\it special} if in convenient holomorphic coordinates,
$$Q(z)=\sum_{j=1}^{n-1} (z_j\overline z_j+\lambda_j Re \ z_j^2), \ \
,\lambda_j\geq 0 \leqno (2)$$

Let $z_j=x_j+iy_j, \ x_j, y_j$ real, $j=1,\ldots,n-1$, then:\vskip 1mm

\noindent(3)\hskip 5mm
$Q(z)=\sum_{l=1}^{n-1}
\big((1+\lambda_l)x_l^2+(1-\lambda_l)y_l^2\big)+ O(|z|^3)$.\hskip 25mm \vskip 1mm

A flat point $p\in S$ is said to be {\it special elliptic} if
$0\leq\lambda_j<1$ for any $j$.

A flat point
$p\in S$ is said to be {\it special 
k-hyperbolic} if $1<\lambda_j$ for $j\in J\subset \{1,\ldots,n-1\}$ and \
$0\leq\lambda_j<1$ \ for $j\in \{1,\ldots,n-1\}\setminus J\not
=\emptyset$, where
$k$  denotes  
the number of elements of $J$.\vskip 1mm

{\it Special elliptic (resp. $k$-hyperbolic) points are elliptic 
(resp. hyperbolic).} \vskip 1mm

\noindent 6.5. {\it Special hyperbolic points}.
\vskip 1mm

\noindent 6.5.1.
We will not consider {\it special parabolic points} (one $\lambda_j=1$ at
least) which don't appear generically.\vskip 1mm

$S$ being given by (1), let $S_0$ be the quadric of equation  
$w=Q(z)$.
{\it Suppose that $S_0$ is flat at $0$ and that $0$ is a special $k$-hyperbolic point. 
Then, in a neighborhood of \ $0$, and with the above local coordinates, it is CR and
nowhere minimal outside $0$, and the CR orbits of $S_0$ are the $(2n-3)$-dimensional
submanifolds given by $w=const. \not =0$.}

The section $w=0$ of $S_0$ is a real quadratic cone $\Sigma'_0$ in
${\bf R}^{2n}$ whose vertex is 0 and, outside 0, it is a CR orbit $\Sigma_0$ in the neighborhood of 0. \vskip 1mm

\noindent 6.6. {\it Foliation by} CR{\it -orbits in the neighborhood of a
special $1$-hyperbolic point}.\vskip 1mm

We mimic the begining of the proof of 2.4.2. in ([DTZ 05],
[DTZ 09]).\vskip 1mm

\noindent 6.6.1. {\it Local {\rm 2}-codimensional submanifolds}. 

In $\C^3$,
consider the 4-dimensional submanifold
$S$ locally defined by the equation 
$$w=\varphi(z)=Q(z)+O(\vert z\vert^3)\leqno (1)$$
and the 4-dimensional submanifold $S_0$ of equation 
$$w=Q(z)\leqno (4)$$ with
$$Q=(\lambda_1+1)x_1^2-(\lambda_1-1)y_1^2+(1+\lambda_2)x_2^2+(1-\lambda_2)
y_2^2$$
 having a special 1-hyperbolic point at 0, $(\lambda_1>1, 0\leq \lambda_2<1)$,
 and the cone $\Sigma'_0$ whose equation is: $Q=0$.
On $S_0$, a CR orbit is the 3-dimensional submanifold ${\cal K}_{w_0}$
whose equation is $w_0=Q(z)$. If $w_0>0$, ${\cal K}_{w_0}$ does not cut
the line $L=\{x_1=x_2=y_2=0\}$; if $w_0<0$, ${\cal K}_{w_0}$ cuts $L$ at
two points.

\vskip 1mm
\noindent 6.6.2. {\it Remark}. $\Sigma_0=\Sigma'_0\setminus 0$ {\it has two
connected components in a neighborhood of} 0.
\vskip 1mm

\noindent {\it Proof}. The equation of  $\Sigma'_0\cap\{y_1=0\}$ is 

$(\lambda_1+1)x_1^2+(1+\lambda_2)x_2^2+(1-\lambda_2)y_2^2=0$ \
\ whose only zero , in the neighborhood of $0$, is $\{ 0\}$: the connected
components are obtained for $y_1>0$ and
$y_1<0$ respectively.   \hskip  25mm $\square$

\vskip 1mm

\noindent 6.6.3. {\it Behaviour of local CR orbits.} 

Under the notations of [DTZ 09], follow the construction of the complex
tangent space $E(z,\varphi(z))$ to the CR orbit at $z$; compare with
$E_0(z,Q(z))$. We know the integral manifold, the orbit  of $E_0(z,Q(z))$;
deduce an evaluation  of the integral manifold of
$E(z,\varphi(z))$. 

\noindent 6.6.4. {\smc Lemma}. {\it Under the above hypotheses, if $k=1$,
the local orbit $\Sigma$ corresponding to $\Sigma_0$} 
 {\it has two connected components in the neighborhood of $0$.}  \vskip 1mm
 
 \noindent {\it Proof}. Use Remark 6.6.2 and the adaptation of the
technique of [DTZ 09].\hskip 5cm  $\square$\vskip 1mm
\vskip 1mm

\noindent 6.7. CR{\it -orbits near a subvariety containing a special
$1$-hyperbolic point}.\vskip 1mm

  \noindent 6.7.2. {\smc Proposition}. {\it Assume that $S\subset\C^n$
$(n\geq 3)$, is a locally closed $(2n-2)$-submanifold, nowhere minimal at
all its CR points, which has a unique spcial $1$-hyperbolic flat complex
point $p$, and such that: 

$(i)$ the orbit $\Sigma$ whose closure $\Sigma'$ contains $p$ is compact;

$(ii)$ $\Sigma$ has two connected components $\sigma_1$,
$\sigma_2$, whose closures are homeomorphic to spheres of dimension $2n-3$.

\vskip 1mm

Then, there exists a neighborhood \ $V$ of \ $\Sigma'$ such that \ $V\setminus
\Sigma'$ is foliated by compact real $(2n-3)$-dimensional CR orbits whose equation,
in a neighborhood of \ $p$ is $(3)$, and, the $w(=x_n)$-axis being assumed to be 
vertical, each orbit being diffeomorphic to 

the sphere ${\bf S}^{2n-3}$ above $\Sigma'$, 

the union of two  spheres ${\bf S}^{2n-3}$ under $\Sigma'$, 

\noindent and there exists a smooth function $\nu$,
 having the CR orbits as the level surfaces. }\

\vskip 1mm

\noindent 6.8. {\it Geometry of the complex points of $S$.}\vskip 1mm

\noindent 6.8.1.  Let $G$ be the manifold of the oriented real linear 
$(2n-2)$-subspaces of $\C^n$. The submanifold $S$ of $\C^n$ has a given orientation
which defines an orientation of the tangent space to $S$ at any point $p\in
S$. By mapping each point of $S$ into its oriented tangent space, we get a smooth
Gauss map
$$t:S\rightarrow G$$

\noindent 6.8..2. {\it Dimension of $G$.} dim $G=2(2n-2)$.\vskip 1mm

\noindent 6.8..3. {\bf Proposition}. {\it For $n\geq 2$, in general, $S$ has
isolated   complex points.} 
\vskip 1mm

\noindent {\it Proof}. Let $\pi\in G$ be a complex hyperplane of $\C^n$ whose 
orientation is induced by its complex structure; the set of such $\pi$ is  
$H=G_{n-1,n}^{\C}=\C{\rm P}^{n-1*}\subset G$, as real submanifold. If $p$ is a
complex point of $S$, then $t(p)\in H$ or $-t(p)\in H$. The set of  complex points of
$S$ is the inverse image by
$t$ of the intersections
$t(S)\cap H$ and $-t(S)\cap H$ in $G$. Since dim $t(S)= 2n-2$, dim $H=2(n-1)$, dim $G=
2(2n-2)$, the intersection  is 0-dimensional, in general.\vskip 1mm

\noindent 6.8.4. {\it Homology of $G$.} (cf [P 08]). $G$ has the structure of a complex quadric;
 let $S_1,S_2$ be generators of $H_{2n-2}(G,\Z)$; we assume
that $S_1$ and
$S_2$ are fundamental cycles of complex projective subspaces of complex dimension $(n-1)$
of $G$. Then, denoting also $S$, the fundamental cycle
of the submanifold $S$ and $t_*$ the homomorphism defined by $t$, we have:
$$t_*(S)\sim u_1S_1+u_2S_2$$where $\sim$ means {\it homologous to}.

 \noindent 6.8.5. {\bf Lemma} (proved for $n=2$ in [CS 51]). {\it With the
notations of
\ $6.8.1'$\ , we have:
\ \ $u_1=u_2$;
\ \ $u_1+u_2=\chi(S)$, Euler-Poincar\'e characteristic of $S$.}

The proof for $n=2$ works for any $n\geq 3$.\vskip 1mm

\noindent 6.8.6. {\it Local intersection numbers of $H$ and $t(S)$ when all
complex points are flat.}

 {\bf Proposition} (known for $n=2$ [Bi 65], {\it here for
$n\geq 3$}). {\it   Let $S$ be a smooth, oriented, compact, 2-codimensional, real
submanifold of $\C^n$ whose all complex points are flat and special. Then, on $S$,  $\sharp$
(special elliptic points)
+ 
\ 
$\sharp$  (special $k$-hyperbolic points, with $k$ even) - $\sharp$ 
(special $k$-hyperbolic points, with $k$ odd) $=\chi(S)$. If $S$ is a 
sphere, this number is 2.}\vskip 2mm

\noindent {\bf 7. Levi-flat hypersurfaces with prescribed boundary: 
particular cases.}\vskip 2mm

\noindent 7.1. To solve the boundary problem by Levi-fllat hypersurfaces, $S$ has to satisfy necessary
and  sufficient local conditions. A way to prove that these conditions can
occur is to construct an example for which  the solution is obvious.
\vskip 1mm
\noindent 7.2. {\it Sphere with elliptic points.}
\vskip 1mm

\noindent 7.2.1. {\it Example}. In $\C^3$, Let $S$ be defined by the
equations:
$$\left\{\matrix{
     z_1\overline z_1+z_2\overline z_2+z_3\overline z_3 &=& 1\cr
z_3&=& \overline z_3\cr
}\right. \leqno (S)$$\vskip 1mm

We have CR-$dim \ S=1$ except at the points $z_1=z_2=0; z_3=\pm 1$ where CR-$dim \ S=2$. $S$ is the unit
 sphere in $\C^2\times\R$; it bounds the unit ball $M$ in $\C^2\times\R$, which is foliated by the complex balls $\C^2\times\{x_3\}\cap M$. The leaves are relatively compact of real
dimension 4 and are bounded by compact leaves (3-spheres) of a foliation of $M$.
\vskip 1mm

\noindent 
7.2.2. {\bf Theorem} [DTZ 05]. {\it Let $S\subset \C^n$,
$n\geq 3$, be a
compact connected  smooth real $2$-codimensional submanifold
satisfying the conditions 

$(i)$ $S$ is nonminimal at every CR point;

$(ii)$every complex point of $S$ is flat and elliptic and there exists at
least one
such point;

$(iii)$ $S$ does not contain complex manifold of dimension
$(n-2)$.

Then $S$ is a topological sphere, and there exists a Levi-flat
$(2n-1)$-subvariety
$\tilde M\subset\C\times\C^n$ with  boundary $\tilde S$
(in the sense of currents)
such that the natural projection $\pi: \C\times\C^n\to \C^n$
restricts to a bijection which is a CR diffeomorphism between $\tilde S$ and
$S$ outside the complex points of $S$.}\vskip 1mm

\noindent 
7.3. {\it Sphere with one special 1-hyperbolic point (sphere with two
horns).}

\noindent 
7.3.1. {\it Example}.
  In $\C^3$, let $(z_j)$, $j=1,2,3$, be the complex 
coordinates and $z_j=x_j+iy_j$. In ${\bf R}^6\cong\C^3$, consider the
4-dimensional subvariety (with negligible singularities) $S$ defined
by:
$y_3=0$
  
$0\leq x_3\leq 1$; \hskip 2mm
$x_3(x_1^2+y_1^2+x_2^2+y_2^2+x_3^2 -1)+(1-x_3)(x_1^4+y_1^4+x_2^4+y_2^4+4x_1^2
-2y_1^2+x_2^2+y_2^2)=0$\vskip 1mm

$-1\leq x_3\leq 0$; \hskip 2mm$x_3=x_1^4+y_1^4+x_2^4+y_2^4+4x_1^2-2y_1^2+x_2^2+y_2^2$\vskip 1mm
The singular set of $S$ is the 3-dimensional section $x_3=0$ along which the 
tangent space is not everywhere (uniquely) defined.

$S$ being in the real hyperplane $\{y_3=0\}$, the complex tangent spaces to 
$S$ are $\{x_3=x^0\}$ for convenient $x^0$.

The set $S$ will be smoothed along the complement of $0$ (origin of $\C^3$) in
its section by the hyperplane 
$\{x_3=0\}$ by a small deformation leaving $h$ unchanged. {\it In the
following $S$ will denote this smooth submanifold}. \vskip 1mm

From elementary analytic geometry, complex points of $S$ are defined by their coordinates:

\noindent $e_3$: $x_j=0,y_j=0$, $(j=1,2)$, $x_3=1$.

\noindent $h$: $x_j=0,y_j=0$, $(j=1, 2)$, $x_3=0$;

\noindent $e_1,e_2$: $x_1=0, y_1=\pm 1, x_2=0, y_2=0, x_3=-1$.\vskip 1mm

{\smc Lemma}. {\it The complex points are flat and special.
The points $e_1,e_2,e_3$ are special elliptic;
the point $h$ is special 1-hyperbolic.}\vskip 1mm

Remark that the numbers of special elliptic and special hyperbolic points
satisfy the conclusion of Proposition 6.8.6.

\vskip 1mm

\noindent 7.3.1'. {\it Shape of $\Sigma'=S\cap\{x_3=0\}$ in the neighborhood of
the origin 
$0$ of $\C^3$}.\vskip 1mm

\noindent {\smc Lemma}. {\it Under the above hypotheses and notations,

$(i)$  $\Sigma=\Sigma'\setminus 0$ has two connected components $\sigma_1$,
$\sigma_2$.

$(ii)$ The closures of the three connected components of $S\setminus\Sigma'$ 
are submanifolds with boundaries and corners.}\vskip 1mm

\noindent {\it Proof}. $(i)$ The only singular point of $\Sigma'$ is 0. We
work in the ball $B(0,A)$ of $\C^2$ $(x_1, y_1,x_2, y_2)$ for small $A$ and
in the 3-space $\pi_\lambda = \{y_2=\lambda x_2\}$, $\lambda\in \R$. For
$\lambda$ fixed, $\pi_\lambda\cong \R^3(x_1, y_1,x_2)$, and
$\Sigma'\cap\pi_\lambda$ is the cone of equation
$4x_1^2-2y_1^2+(1+\lambda^2)x_2^2+O(|z|^3)=0$ with vertex 0 and basis  in the plane
$x_2=x_2^0$ the hyperboloid
$H_\lambda$ of equation $4x_1^2-2y_1^2+(1+\lambda^2)x_2^{0 2}+O(|z|^3)=0$; the curves
$H_\lambda$ have no common point outside 0. So, when $\lambda$ varies, the surfaces
$\Sigma'\cap\pi_\lambda$ are disjoint outside 0. The set $\Sigma'$ is
clearly connected; $\Sigma'\cap\{y_1=0\}=\{0\}$, the origin of $\C^3$;
from above: $\sigma_1= \Sigma\cap\{y_1>0\}$; $\sigma_2=
\Sigma\cap\{y_1<0\}$. 

$(ii)$ The three connected components of $S\setminus\Sigma'$ are the
components which contain, respectively $e_1$, $e_2$, $e_3$ and whose
boundaries are $\overline \sigma_1$, $\overline \sigma_2$, $\overline
\sigma_1\cup \overline \sigma_2$; these boundaries have corners as shown
in the first part of the proof.\hskip 14cm $\square$\vskip 1mm

The connected component of $\C^2\times \R\setminus S$
containing the point $(0,0,0,0,1/2) $ is the Levi-flat solution, the complex
leaves being the sections by the hyperplanes $x_3=x_3^0$, $-1<x_3^0<1$.

The sections by the hyperplanes $x_3=x_3^0$ are diffeomorphic to a 3-sphere for 
$0<x_3^0<1$ and to the union of two disjoint 3-spheres for $-1<x_3^0<0$, as can be
shown intersecting $S$ by lines through the origin in the
hyperplane $x_3=x_3^0$;
$\Sigma'$ is homeomorphic to the union of two 3-spheres with a common point.
 \vskip 1mm

\noindent 
7.3.2. {\bf Proposition} (cf [D 08],  Proposition 2.6.1). {\it Let
$S\subset \C^n$ be a compact connected real 2-codimensional manifold such that
the following holds:

$(i)$ $S$ is a topological sphere; $S$ is nonminimal at every CR point;

$(ii)$ every complex point of $S$ is flat; there exist three special elliptic
points $e_j, j=1,2,3$
and one special $1$-hyperbolic point $h$;

$(iii)$ $S$ does not contain complex manifolds of dimension
$(n-2)$;

$(iv)$  the singular {\rm CR} orbit $\Sigma'$ through $h$ on $S$ is compact and
$\Sigma'\setminus \{h\}$ has two connected components $\sigma_1$ and
$\sigma_2$ whose closures are homeomorphic to spheres of dimension $2n-3$;

$(v)$ the closures $S_1, S_2, S_3$ of the three connected components $S'_1,
S'_2, S'_3$ of
$S\setminus\Sigma'$  are submanifolds with (singular) boundary.

Then each $S_j\setminus \{e_j\cup\Sigma'\}$, $j=1,2,3$
carries a foliation ${\cal F}_j$ of class
$C^\infty$ with $1$-codimensional {\rm CR} orbits as compact leaves.}
\vskip 1mm

{\it Proof}. From conditions (i) and (ii), $S$ satisfying
the hypotheses of Proposition 6.3.2, near any elliptic
flat point $e_j$, and of Proposition 6.7.2 near $\Sigma'$, all CR orbits are
diffeomorphic to the sphere ${\bf S}^{2n-3}$.
The assumption (iii) guarantees that
all CR orbits in $S$ must be of real dimension $2n-3$.
Hence, by removing small connected open saturated neighborhoods
of all special elliptic points, and of $\Sigma'$,
we obtain, from $S\setminus\Sigma'$, three compact manifolds $S_j"$,
$j=1,2,3$,
 with boundary and 
with the foliation ${\cal F}_j$ of codimension $1$ given by its CR
orbits whose first
cohomology group with values in ${\bf R}$ is 0, near $e_j$.
It is easy to show that this foliation is transversely oriented.

\noindent 7.3.2'. Recall the Thurston's Stability Theorem 
($\lbrack$ CaC$\rbrack$, Theorem 6.2.1).
{\it Let $(M,{\cal F})$ be a compact, connected, transversely-orientable,
foliated manifold with boundary or corners, of codimension 1, of class
$C^1$. If there is a compact leaf $L$ with $H^1(L,{\bf R})=0$, then every leaf
is homeomorphic to $L$ and $M$ is homeomorphic to $L\times \lbrack
0,1\rbrack$, foliated as a product}, 

Then, from the above theorem,  $S_j"$ is homeomorphic to
${\bf S}^{2n-3}\times [0,1]$ with CR orbits being of the form
${\bf S}^{2n-3}\times\{x\}$ for $x\in [0,1]$. Then the full manifold $S_j$ is
homeomorphic to a half-sphere supported by ${\bf S}^{2n-2}$ and
${\cal F}_j$ extends to $S_j$; $S_3$ having its boundary pinched at the point $h$.\vskip 1mm

\noindent 
7.3.3.  {\bf Theorem.}
{\it Let $S\subset \C^n$, $n\geq 3$, be a
compact connected  smooth real $2$-codimensional submanifold
satisfying the conditions $(i)$ to $(v)$ of Proposition 7.3.2. 
Then there exists a Levi-flat $(2n-1)$-subvariety $\tilde M\subset\C\times\C^n$
with  boundary $\tilde S$
(in the sense of currents)
such that the natural projection $\pi: \C\times\C^n\to \C^n$
restricts to a bijection which is a CR diffeomorphism between $\tilde S$ and
$S$ outside the complex points of $S$.}

{\it Proof}. 
By Proposition 6.3.2 , for every $e_j$, a
continuous function
$\nu'_j$, $C^{\infty}$ outside $e_j$, can be constructed in a
neighborhood $U_j$ of $e_j$, $j=1,2,3$, and by Proposition 6.7.2, we have an
analogous result   in a neighborhood of $\Sigma'$.

Furthermore, from section 7.3.2', a smooth function $\nu"_j$ whose level
sets are the leaves of ${\cal F}_j$ can be obtained
globally on
$S'_j\setminus  \{e_j\cup\Sigma'\}$. With the functions $\nu'_j$ and $\nu"_j$, and
analogous functions near $\Sigma'$, then using a partition of unity, we obtain a
global smooth function
$\nu_j\colon S_j\to {\bf R}$ without critical points away from the complex points
$e_j$ and from $\Sigma'$.

Let $\sigma_1$, resp. $\sigma_2$ the two connected, relatively compact 
components of $\Sigma\setminus \{h\}$, according to condition $(iv)$;
$\overline\sigma_1$, resp. $\overline\sigma_2$ are the boundary of $S_1$,
resp. $S_2$, and $\overline \sigma_1\cup\overline\sigma_2$ the boundary of
$S_3$. We can assume that the three functions $\nu_j$ are finite valued and get the same values on $\overline\sigma_1$ and $\overline\sigma_2$.
Hence a function $\nu: S\rightarrow{\bf R}$.

The submanifold $S$ being, locally, a boundary of a Levi-flat hypersurface, is
orientable. We now set $\tilde S=N={\rm gr}\,\nu = \{(\nu(z),z) : z\in S\}$.
Let $S_s=\{e_1,e_2, e_3,\overline{\sigma_1\cup\sigma_2}\}$.

$\lambda: S\rightarrow\tilde S\hskip 2mm\big( z\mapsto \nu ((z),z)\big)$ is
bicontinuous; $\lambda\vert_{S\setminus S_s}$ is a diffeomorphism; moreover
$\lambda$ is a CR map. Choose an orientation on $S$. Then $N$ is an (oriented)
CR subvariety with the negligible set of singularities $\tau=\lambda(S_s)$.

At every point of $S\setminus S_s$, $d_{x_1} \nu\not = 0$, then
condition (H) (section 5.1.1) is satisfied at every point of
$N\setminus\tau$.

Then all the assumptions of Theorem 5.1.2 being satisfied by
$N=\tilde S$, in a particular case, we conclude that $N$ is the
boundary of a Levi-flat
$(2n-2)$-variety (with negligible singularities) $\tilde M$ in ${\bf  
R}\times\C^n$.\vskip 1mm

Taking $\pi: \C\times \C^n\to \C^n$ to be the standard projection,
we obtain the conclusion. \vskip 1mm

\noindent 
7.4. {\it  Case of a torus.}
\vskip 1mm

\noindent 
7.4.1. {\it Euler-Poincar\'e characteristic of a torus is $\chi({\bf
T}^k)=0$.}
\vskip 1mm

\noindent 
7.4.2. {\it Example}. In $ÖC^3$, let $(z_j)$, $j=1,2,3$, be the 
complex 
coordinates and $z_j=x_j+iy_j$. In ${\bf R}^6\cong\C^3$, consider the
4-dimensional subvariety (with negligible singularities) $S$ defined
by:

$y_3=0$

$0\leq x_3\leq 1$; \hskip 2mm
$x_3(x_1^2+y_1^2+x_2^2+y_2^2+x_3^2 -1)+(1-x_3)(x_1^4+y_1^4+x_2^4+y_2^4+4x_1^2
-2y_1^2+x_2^2+y_2^2)=0$

$-\frac{1}{2} \leq x_3\leq 0$; \hskip
2mm$x_3=x_1^4+y_1^4+x_2^4+y_2^4+4x_1^2-2y_1^2+x_2^2+y_2^2$

\noindent glue it with the symmetric with respect to the real hyperplane
$x_3=-\frac{1}{2}$, and {\it and smooth along} 
$\{x_3=0\}$, $\{x_3=\pm\frac{1}{2}\}$. The complex points are flat and
special.\vskip 1mm

\noindent 
7.4.3. {\bf Theorem.} {\it Let $S\subset \C^n$, $n\geq 3$, be a
compact connected  smooth real $2$-codimensional submanifold
satisfying the following conditions:

$(i)$ $S$ is a topological torus; $S$ is nonminimal at every CR point;

$(ii)$ every complex point of $S$ is flat; there exist two special elliptic points $e_1,  e_2$ and two special $1$-hyperbolic points $h_1, h_2$;

$(iii)$ $S$ does not contain complex manifolds of dimension $(n-2)$;

$(iv)$  the singular \ {\rm CR} orbits $\Sigma'_1, \Sigma'_2$ through
$h_1$ and $h_2$ on $S$ are compact and, for $j=1,2$, $\Sigma'_j\setminus \{h_j\}$ have two connected components $\sigma^1_j$ and $\sigma^2_j$;

$(v)$ the closures $S_1, S_2, S_3, S_4$ of the four connected components
$S'_1, S'_2, S'_3, S'_4$ of $S\setminus\Sigma'_1\cup\Sigma'_2$  \ are submanifolds with (singular) boundary.

Then there exists a Levi-flat $(2n-1)$-subvariety $\tilde M\subset\C\times\C^n$
with  boundary $\tilde S$ (in the sense of currents)
such that the natural projection $\pi: \C\times\C^n\to \C^n$
restricts to a bijection which is a CR diffeomorphism between $\tilde S$ and
$S$ outside the complex points.}\vskip 1mm

\noindent 
7.5. {\it Generalizations.}\vskip 1mm

\noindent  7.5.1. {\it Elementary models and their gluing.}
 The examples and the proofs of the theorems 
when $S$ is homeomorphic to a sphere (sections 7.3) or a torus (section  7.4)
suggest the following definitions.\vskip 1mm

\noindent 7.5.2. {\it Definitions}.\vskip 1mm

Let $T'$ be a smooth, locally closed (i.e. closed in an open set), connected submanifold of $\C^n$, $n \geq 3$. We assume
that
$T'$ has the following properties:

$(i)$  $T'$ is relatively compact, non necessarily compact, and of codimension 2. 

$(ii)$ $T'$ is nonminimal at every CR point;

$(iii)$ $T'$ has exactly 2 complex points which are flat and either special elliptic or special 1-hyperbolic.

$(iv)$ If $p\in T'$ is 1-hyperbolic, the singular orbit $\Sigma'$
through
$p$ is compact, $\Sigma'\setminus p$ has two connected components
$\sigma_1$,
$\sigma_2$, whose closures are homeomorphic to spheres of dimension $2n-3$.

$(v)$  If $p\in T'$ is 1-hyperbolic, in the neighborhood of
$p$, with convenient coordinates, the equation of $T'$, up to third order terms
is 
$$z_n=\sum_{j=1}^{n-1} (z_j\overline z_j+\lambda_j {\cal R}e \ z_j^2); \ 
\lambda_1>1; \ 0\leq\lambda_j<1 \ \ {\rm for} \ \  j\not =1 $$

or in real coordinates $x_j,y_j$
 with $z_j=x_j+iy_j$, \quad
$$x_n=\big((\lambda_1+1)x_1^2-(\lambda_1-1)y_1^2\big)+\sum_{j=2}^{n-1}
\big((1+\lambda_j)x_j^2+(1-\lambda_j)y_j^2\big)+ O(|z|^3)$$
Other configurations are easily imagined.\vskip 1mm

\noindent {\it up- and down-} 1{\it -hyperbolic points}. Let $T$ be
the $(2n-2)$-submanifold with (singular) boundary contained into $T'$ 
such that either $\overline\sigma_1$ (resp. $\overline\sigma_2$) is the
boundary of $T$ near $p$, or $\Sigma'$ is the boundary of $T$ near $p$.
In the first case, we say that $p$ is 1-{\it up}, (resp. 2-{\it up}), in
the second that $p$ is {\it down}. Such a $T$ will be called an
{\it elementary model}.\vskip 1mm

For instance, $T$ is 1-up and has one special elliptic point, we solve the boundary problem as in $S_1$ in the proof of Theorem 7.3.3.\vskip 1mm

\noindent 7.5.3.  The gluing (to be precised) happens between two compatible
elementary models along boundaries, for instance down and 1-up.\vskip 1mm

\noindent 7.6. {\it Other possible generalizations.}\vskip 1mm
The mixed Plateau problem can be set up in projective space $\cp n$ and in
subspaces of $\cp n$ on which the complex Plateau problem can be solved, using Statement 5.2.6, its gemeralisation to any $n\geq 3$ and a better geometric condition on the given boundary.

\vskip 5mm

\centerline {\bf References}\vskip 2mm

 [BeK 91] E. Bedford \& W. Klingenberg, {\it On the
envelopes of
holomorphy of a 2-sphere in $\C^2$}, J. Amer. Math. Soc. {\bf 4}
(1991), 623-646.

 [Bi 65] E. Bishop, {\it Differentiable manifolds in complex Euclidean
space}, Duke Math. J. {\bf 32} (1965), 1-22.

[CaC]  A. Candel \& L. Conlon, {\it Foliations. I.}
Graduate Studies in Mathematics, {\bf 23}. American Mathematical Society,
Providence, RI, 2000.

[CS 51] S.S. Chern and E. Spanier, {\it A theorem on orientable
surfaces in four-dimensional space}, Com. Math. Helv., {\bf 25} (1951),
205-209.

$\lbrack$Di 98$\rbrack$ \hskip 2mm T.C. Dinh, Enveloppe polynomiale d'un compact de
longueur finie
 et cha\^{\i}nes holomorphes \`a bord rectifiable, Acta Math. {\bf 180} (1998), 31-67.

[D 08] P. Dolbeault, {\it On Levi-flat hypersurfaces with given
boundary in $\C^n$}, Science in China, Series A: Mathematics, 
Apr. 2008, Vol.
51, No 4, 551-562. 

[D 09] P. Dolbeault, {\it On Levi-flat hypersurfaces with given
boundary: special hyperbolic points}, in preparation. 
 
[DTZ 05] P. Dolbeault, G. Tomassini, D. Zaitsev, {\it On
boundaries of Levi-flat hypersurfaces in $\C^n$}, C.R. Acad. Sci. Paris,
Ser. I {\bf 341} (2005) 343-348. 

[DTZ 09] P. Dolbeault, G. Tomassini, D. Zaitsev, {\it On Levi-flat 
hypersurfaces with prescribed boundary}, Pure and Applied Math. Quarterly
 {\bf 6}, N. 3\ (Special Issue: In honor  of \ Joseph J. Kohn), 725---753, (2010) arXiv:0904.0481
 
$\lbrack$DH 97$\rbrack$ P. Dolbeault et G. Henkin, {\it Chaines
holomorphes de bord donn\'e dans $\C P^n$}, Bull. Soc. Math. France, {\bf 125}, 383-445.

[H 77] R. Harvey, {\it Holomorphic chains and their boundaries}, Proc. Symp.
Pure Math. {\bf 30}, Part I, Amer. Math. Soc. (1977),309-382.

$\lbrack$HL $75 \rbrack$ \hskip 2mm R. Harvey and B. Lawson, On boundaries of complex analytic
varieties, I, Ann. of Math., {\bf 102},  (1975), 233-290.

$\lbrack$HL $77 \rbrack$ \hskip 2mm R. Harvey and B. Lawson, On boundaries of complex analytic
varieties, II, Ann. of Math., {\bf 106}, (1977), 213-238.

[HL 04] Harvey, F. Reese; Lawson, H. Blaine, Jr. Boundaries of varieties in projective manifolds. J. Geom. Anal. 14 (2004), no. 4, 673--695.

$\lbrack$Ki $79 \rbrack$ \hskip 2mm J. King, Open problems in geometric function theory, Proceedings of the fifth international 
 symposium of Math.,p. 4, The Taniguchi foundation, 1978.

$\lbrack$Lce 95$\rbrack$\hskip 2mm M.G. Lawrence, Polynomial hulls of rectifiable curves, Amer. J. Math., {\bf 117} (1995), 405-417. 

[P 08] P. Polo, Grassmanniennes orient\'ees r\'eelles, e-mail personnelle, 21
f\'ev. 2008 

$\lbrack$Rs 59$\rbrack$\hskip 2mm   W. Rothstein, Bemerkungen zur Theorie komplexer R\"aume, Math.
Ann., {\bf 137} (1959), 304-315.

$\lbrack$We 58$\rbrack$\hskip 2mm J. Wermer, The hull of a curve in $\C^n$. Ann.of Math. {\bf
68}  (1958), 550-561.

\end